\documentclass[11pt,a4paper,fleqn]{article}
\usepackage{a4wide,amsfonts,amsmath,latexsym,amssymb,euscript,graphicx,units,mathrsfs}

\usepackage{graphicx}
\usepackage{color}
\usepackage{amssymb}
\usepackage{amssymb}
\usepackage[T1]{fontenc}
\usepackage{latexsym}
\usepackage{xypic}
\usepackage{eufrak}
\usepackage{euscript}
\usepackage{amsfonts,amsmath}
\usepackage{verbatim}
\usepackage{fancyhdr}
\usepackage[english]{babel}
\usepackage{mathrsfs}
\usepackage{units}
\usepackage{hyperref}

\newtheorem{prop}{Proposition}[section]

\newtheorem{thm}[prop]{Theorem}
\newtheorem{defi}[prop]{Definition}

\renewcommand{\geq}{\geqslant}
\def\leq{\leqslant}

\newcommand{\Z}{\mathbb{Z}}
\newcommand{\R}{\mathbb{R}}

\def\HH{\EuFrak H}

\def\1{{\mathbf{1}}}

\def\sk{{\mathbb{D}}}

\def\1{{\mathbf{1}}}
\def\0.5{{\frac{1}{2}}}

\newcommand{\qed}{\nopagebreak\hspace*{\fill}
{\vrule width6pt height6ptdepth0pt}\par}

\begin{document}

\title{\textbf{ 
The optimal fourth moment theorem
}}
\author{Ivan Nourdin and Giovanni Peccati}
\maketitle
\abstract{We compute the exact rates of convergence in total variation associated with the `fourth moment theorem' by Nualart and Peccati (2005), stating that a sequence of random variables living in a fixed Wiener chaos  verifies a central limit theorem (CLT) if and only if the sequence of the corresponding fourth cumulants converges to zero. We also provide an explicit illustration based on the Breuer-Major CLT for Gaussian-subordinated random sequences.\\


}

\section{Introduction and main result}

On a suitable probability space $(\Omega, \mathscr{F}, P)$, let $X = \{X(h) : h\in \HH\}$ be an isonormal Gaussian process over a real separable Hilbert space $\HH$, and let $\{C_q : q=0,1,...\}$ be the sequence of Wiener chaoses associated with $X$ (see Section \ref{ss:mall} for details). The aim of the present work is to prove an optimal version of the following quantitative `fourth moment' central limit theorem (CLT) --- which combines results from  \cite{np-ptrf} and \cite{NunuGio}. 

\begin{thm}[Fourth moment theorem \cite{np-ptrf, NunuGio}]\label{t:fmt} Fix an integer $q\geq 2$ and let $N\sim \mathscr{N}(0,1)$ be a standard Gaussian random variable. Consider a sequence $\{F_n : n\geq 1\}$ living in the $q$th Wiener chaos $C_q$ of $X$, and assume that $E[F_n^2] =1$. Then, $F_n$ converges in distribution to $N\sim \mathscr{N}(0,1)$ if and only if $E[F_n^4]\to 3  = E[N^4]$. Also, the following upper bound on the total variation distance holds for every $n$:
\begin{equation}\label{e:b} 
d_{TV}(F_n,N) \leq \sqrt{\frac{4q-4}{3q} } \sqrt{|E[F_n^4]-3|}. 
\end{equation}
\end{thm}

Recall that, given two real-valued random variables $Y,Z$, the {\it total variation distance} between the law of $Y$ and $Z$ is given by 
\begin{eqnarray}
d_{TV} (Y,Z) &=& \sup_{A\in \mathscr{B}(\R)} \left| P[Y\in A] -P[Z\in A] \right| \label{e:tv1}\\
&=& \frac12\sup \left| E[g(Y)] - E[g(Z)]\right|, \label{e:tv2}
\end{eqnarray}
where, in view e.g. of Lusin's Theorem, the supremum is taken taken over the class of continuously differentiable functions $g$ that are bounded by 1 and have compact support. Also, if the distributions of $Y,Z$ have densities (say $f_Y, f_Z$) then one has the integral representation
\begin{equation}\label{e:tv3}
d_{TV}(Y,Z) = \frac12\int_\R \left| f_Y(t) - f_Z(t)\right| dt.
\end{equation}
It is a well-known fact that the topology induced by $d_{TV}$, over the class of all probabilities on $\R$, is strictly stronger than the topology of convergence in distribution (see e.g. \cite{dudley_book}). It is also important to notice that bounds analogous to \eqref{e:b} hold for other distances, like for instance the Kolmogorov and $1$-Wasserstein distances (see \cite[Chapter 5]{NP11}).

\medskip

The first part of Theorem \ref{t:fmt} was first proved in \cite{NunuGio} by using tools of continuous-time stochastic calculus. The upper bound \eqref{e:b} comes from reference \cite{np-ptrf}, and is obtained by combining the {\it Malliavin calculus of variations} (see e.g. \cite{NP11, Nu06}) with the {\it Stein's method} for normal approximations (see \cite{ChGoSh11, NP11}). The content of Theorem \ref{t:fmt} has sparkled a great amount of generalisations and applications, ranging from density estimates \cite{NV-ejp} to entropic CLTs \cite{NoPeSw}, and from estimates for Gaussian polymers \cite{viensSPA} to universality results \cite{NPRaop}. The reader is referred to the monograph \cite{NP11} for a self-contained introduction to the theoretical aspects of this direction of research, as well as to \cite{MP11, NouBookFBM} for applications, respectively, to the high-frequency analysis of fields defined on homogeneous spaces, and to the power variation of stochastic processes related to fractional Brownian motion.  One can also consult the constantly updated webpage 
\begin{center}
\url{http://www.iecn.u-nancy.fr/~nourdin/steinmalliavin.htm}  
\end{center}
for literally hundreds of results related to Theorem \ref{t:fmt} and its ramifications.

\medskip

One challenging question is whether the upper bound (on the rate of convergence in total variation) provided by \eqref{e:b} is optimal, or rather it can be ameliorated in some specific situations. Recall that a positive sequence $\{ \phi(n) : n\geq 1\}$ decreasing to zero yields an {\it optimal rate of convergence}, with respect to some suitable distance $d(\cdot, \cdot)$, if there exist two finite constants $0<c<C$ (independent of $n$) such that
\[
c\, \phi(n) \leq d(F_n, N) \leq C\, \phi(n),\quad \mbox{ for all }\,  n\geq 1.
\]
There are indeed very few references studying optimal rates of convergence for CLTs on a Gaussian space. Our paper \cite{np-aop} provides some partial characterisation of optimal rates of convergence in the case where $d = d_{Kol}$ (the Kolmogorov distance), whereas \cite{campese} contains similar findings for multidimensional CLTs and distances based on smooth mappings. Of particular interest for the present analysis is the work \cite{BBNP}, that we cowrote with H. Bierm\'e and A. Bonami, where it is proved that, whenever the distance $d$ is defined as the supremum over a class of {\it smooth enough} test functions (e.g., twice differentiable with a bounded second derivative), an optimal rate of convergence in Theorem \ref{t:fmt} is given by the sequence
\begin{equation}\label{e:q}
{\bf M}(F_n) := \max\left\{  |E[F_n^3] | ,\, | E[F_n^4] - 3 | \right\}, \quad n\geq 1.
\end{equation}
In particular, if $E[F_n^3] = 0$ (for instance, if $q$ is odd), then the rate suggested by the estimate \eqref{e:b} is suboptimal by a whole square root factor. We observe that, if $F$ is a non-zero element of the $q$th Wiener chaos of $X$, for some $q\geq 2$, then $E[F^4] > 3E[F^2]^2$; see e.g. \cite{NunuGio} or \cite[Chapter 5]{NP11}.

\medskip

The statement of the subsequent Theorem \ref{t:1}, which is the main achievement of the present paper, provides a definitive characterisation of the optimal rate of convergence in the total variation distance for the fourth moment CLT appearing in Theorem \ref{t:fmt}. 
It was somehow unexpected that our optimal rate only 
relies on the simple quantity (\ref{e:q}) and also that we do not need to impose a 
further restriction than being an element of a given Wiener chaos. 
One remarkable consequence of our findings is that this optimal rate exactly coincides with the one related to the smooth test functions considered in \cite{BBNP}. We will see that the proof relies on estimates taken from the paper \cite{BBNP}, that we combine with a new fine analysis of the main upper bound proved in \cite{np-ptrf} (see Proposition \ref{p:tv}).

\begin{thm}[Optimal fourth moment theorem in total variation]\label{t:1} Fix $q\geq 2$. Let $\{F_n : n\geq 1\}$ be a sequence of random variables living in the $q$th Wiener chaos of $X$, such that $E[F_n^2] =1$. Then, $F_n$ converges in distribution to $N\sim \mathscr{N}(0,1)$ if and only if $E[F_n^4]\to 3  = E[N^4]$. In this case, one has also that $E[F_n^3]\to 0$ and there exist two finite constants $0<c <C$ (possibly depending on $q$ and on the sequence $\{F_n\}$, but not on $n$) such that the following estimate in total variation holds for every $n$:
\begin{equation}\label{e:mainest}
c\,  {\bf M} (F_n) \leq d_{TV}(F_n,N) \leq C\,  {\bf M} (F_n),
\end{equation}
where the quantity ${\bf M} (F_n)$ is defined according to (\ref{e:q}). 
\end{thm}

The rest of the paper is organised as follows. Section \ref{s:pre}
contains some notation and useful preliminaries. Our main result, Theorem \ref{t:1},
is proved in Section 3. Finally,  in Section 4 we provide an explicit illustration based on the Breuer-Major CLT for Gaussian-subordinated random sequences.

\section{Notation and preliminaries}\label{s:pre}

\subsection{Cumulants}

In what follows, the notion of \textit{%
cumulant} is sometimes used. Recall that, given a random variable $Y$ with finite moments of
all orders and with characteristic function $\psi _{Y}\left( t\right) =%
E\left[ \exp \left( itY\right) \right] $ ($t\in \mathbb{R}$), one
defines the sequence of cumulants (sometimes known as \textit{semi-invariants%
}) of $Y$, noted $\left\{ \kappa _{n}\left( Y\right) :n\geq 1\right\} $, as
\begin{equation}
\kappa_{n}\left( Y\right) =\left( -i\right) ^{n}\frac{d^{n}}{dt^{n}}\log \psi
_{Y}\left( t\right) \mid _{t=0}\text{, \ \ }n\geq 1\text{.}  \label{Eq : cum}
\end{equation}%
For instance, $\kappa_{1}\left( Y\right) =E[ Y] $, $\kappa_{2}\left( Y\right) =E\left[ (Y-E[Y])^2\right]={\rm Var}\left( Y\right) $, and, if $E[Y]=0$,
\begin{equation*}
\kappa_{3}\left( Y\right) =E[ Y^{3}]\quad \mbox{and}\quad \kappa_{4}\left( Y\right) =E[ Y^{4}] -3E[Y^2]^2.
\end{equation*}%

\subsection{Hermite polynomials}

\noindent We will also denote by $\{H_q : q=0,1,...\}$ the sequence of {\it Hermite polynomials} given by the recursive relation $H_0=1$ and $H_{q+1}(x) = xH_q(x) - H'_q(x)$, in such a way that $H_1(x) = x$, $H_2(x) = x^2-1$, $H_3(x) = x^3-3x$ and $H_4(x) = x^4 - 6x^2+3$. We recall that Hermite polynomials constitute a complete orthogonal system of the space $L^2\big(\R, \mathscr{B}(\R), (2\pi)^{-1/2} e^{-x^2/2}dx\big)$.

\subsection{Stein's equations}

See \cite[Chapter 3]{NP11} for more details on the content of this section. Let $N\sim \mathscr{N}(0,1)$. Given a bounded and continuous function $g: \R\to \R$, we define the {\it Stein's equation} associated with $g$ to be the ordinary differential equation
\begin{equation}\label{e:steine}
f'(x) -xf(x) = g(x) - E[g(N)].
\end{equation}
It is easily checked that every solution to \eqref{e:steine} has the form $ce^{x^2/2}+ f_g(x)$, where 
\begin{equation}\label{e:j}
f_g(x) = e^{x^2/2} \int_{-\infty}^x \{g(y) - E[g(N)]\} e^{-y^2/2}dy, \quad x\in \R.
\end{equation}
Some relevant properties of $f_g$ appear in the next statement (see e.g. \cite[Section 3.3]{NP11} for proofs).

\begin{prop}[Stein's bounds]\label{p:stein} Assume $g:\R\to\R$ is continuous and bounded. Then $f_g$ given by (\ref{e:j}) is $C^1$ and satisfies
\[ 
\|f_g\|_\infty \leq \sqrt{\pi/2} \, \| g - E[g(N)]\|_\infty\quad \mbox{and} \quad \|f'_g\|_\infty \leq 2\, \| g - E[g(N)]\|_\infty.
\] 
\end{prop}

\subsection{The language of Gaussian analysis and Malliavin calculus}\label{ss:mall}

We now briefly recall some basic notation and results connected to Gaussian analysis and Malliavin calculus. The reader is referred to \cite{NP11, Nu06} for details.

\medskip

Let $\HH$ be a real separable Hilbert space with inner product $\langle \cdot, \cdot \rangle_\HH$. Recall that an {\it isonormal Gaussian process} over $\HH$ is a centered Gaussian family $X = \{X(h) : h\in \HH\}$, defined on an adequate probability space $(\Omega, \mathscr{F}, P)$ and such that $E[X(h) X(h')] = \langle h, h'\rangle_\HH$, for every $h,h'\in \HH$. For every $q=0,1,2,...$, we denote by $C_q$ the $q$th {\it Wiener chaos} of $X$. We recall that $C_0= \R$ and, for $q\geq 1$, $C_q$ is the $L^2$-closed space composed of those random variables having the form $I_q(f)$ where $I_q$ indicates a {\it multiple Wiener-It\^o} integral of order $q$ and $f\in \HH^{\odot q}$ (the $q$th symmetric tensor power of $\HH$). Recall that $L^2(\sigma(X),P):= L^2(P) = \bigoplus_{q=0}^\infty C_q$, that is: every square-integrable random variable $F$ that is measurable with respect to $\sigma(X)$ (the $\sigma$-field generated by $X$) admits a unique decomposition of the type
\begin{equation}\label{e:hh}
F = E[F] +\sum_{q=1}^\infty I_q(f_q),
\end{equation}
where the series converges in $L^2(P)$, and $f_q \in \HH^{\odot q}$, for $q\geq 1$. This last result is known as the {\it Wiener-It\^o chaotic decomposition} of $L^2(P)$. When the kernels $f_q$ in \eqref{e:hh} all equal zero except for a finite number, we say that $F$ has a {\it finite chaotic expansion}. According to a classical result   (discussed e.g. in \cite[Section 2.10]{NP11}) the distribution of non-zero random variables with a finite chaotic expansion has necessarily a density with respect to the Lebesgue measure.

\medskip

We will use some standard operators from Malliavin calculus. The {\it Malliavin derivative} $D$ has domain $\mathbb{D}^{1,2}\subset L^2(P)$, and takes values in the space $L^2(P;\HH)$ of square-integrable $\HH$-valued random variables that are measurable with respect to $\sigma(X)$. The {\it divergence operator} $\delta$ is defined as the adjoint of $D$. In particular, denoting by ${\rm dom}\, \delta$ the domain of $\delta$, one has the so-called {\it integration by parts formula}: for every $D\in \mathbb{D}^{1,2}$ and every $u\in \delta(u)$,
\begin{equation}\label{e:ipp}
E[F\delta(u)] = E[\langle DF ,u\rangle_\HH].
\end{equation}
We will also need the so-called {\it generator of the Ornstein-Uhlenbeck semigroup}, written $L$, which is defined by the relation
\[
L = - \delta D, 
\] 
meaning that $F$ is in ${\rm dom}\, L$ (the domain of $L$) if and only if $F\in \mathbb{D}^{1,2}$ and $DF\in {\rm dom} \, \delta$, and in this case $LF = -\delta DF$. The {\it pseudo-inverse} of $L$ is denoted by $L^{-1}$. It is important to note that $L$ and $L^{-1}$ are completely determined by the following relations, valid for every $c\in \R$, every $q\geq 1$ and every $f\in \HH^{\odot q}$:
\begin{equation*}
L\, c = L^{-1}\, c = 0, \quad LI_q(f) = -qI_q(f), \quad L^{-1}I_q(f) = -\frac1q I_q(f).
\end{equation*}
We will exploit the following {\it chain rule} (see e.g. \cite[Section 2.3]{NP11}): for every $F\in \mathbb{D}^{1,2}$ and every mapping $\varphi : \R\to \R$ which is continuously differentiable and with a bounded derivative, one has that $\varphi(F)\in \mathbb{D}^{1,2}$, and moreover
\begin{equation}
D\varphi(F) = \varphi'(F)DF.
\end{equation}

\medskip

We conclude this section by recalling the definition of the operators $\Gamma_j$, as first defined in \cite{np-jfa}; see \cite[Chapter 8]{NP11} for a discussion of recent developments.

\begin{defi}[Gamma operators]\label{d:gamma} {\rm Let $F$ be a random variable having a finite chaotic expansion. The sequence of random variables $\{\Gamma_j(F) : j\geq 0\}$ is recursively defined as follows. Set $\Gamma_0(F) = F$
and, for every $j\geq 1$, 
\[
\Gamma_{j}(F) = \langle DF,-DL^{-1}\Gamma_{j-1}(F)\rangle_{\HH}.
\]
In view of the product formulae for multiple integrals (see e.g. \cite[Theorem 2.7.10]{NP11}), each $\Gamma_j(F)$ is a well-defined random variable having itself a finite chaotic expansion.}
\end{defi}

\subsection{Some useful bounds}

The next bound contains the main result of \cite{np-ptrf}, in a slightly more general form (whose proof can be found in \cite{nourdinln}) not requiring that the random variable $F$ has a density.

\begin{prop}[General total variation bound]\label{p:tv} Let $F$ be a centered element of $\sk^{1,2}$, and let $N\sim \mathscr{N}(0,1)$. Then, 
\begin{equation}\label{e:ivangio1}
d_{TV}(F,N) \leq 2E\left| E[1- \langle DF, -DL^{-1}F\rangle_\HH | F] \right|\leq 2E\left| 1- \langle DF, -DL^{-1}F\rangle_\HH\right|. 
\end{equation}
In particular, if $F=I_q(f)$ belongs to the $q$th Wiener chaos $C_q$ and $F$ has unit variance, one has that $\langle DF, -DL^{-1}F\rangle_\HH = \frac{1}{q}\|DF\|^2_\HH$, $E[F^4]>3$, and the following inequality holds:
\begin{equation}\label{e:ivangio2}
E\left| 1- \frac{1}{q}\|DF\|^2_\HH \right| \leq \sqrt{\frac{q-1}{3q} (E[F^4]-3) }.
\end{equation}

\end{prop}

\medskip

\noindent Of course, relation \eqref{e:ivangio2} is equivalent to \eqref{e:b}. In particular, using the language of cumulants, the estimates in the previous statement yield that, for $F=I_q(f)$ ($q\geq 2$) with unit variance,
\begin{equation}\label{NPfirst}
d_{TV}(F,N)\leq \frac{2}{\sqrt{3}}\sqrt{\kappa_4(F)}.
\end{equation}

The following result (which is taken from \cite{BBNP}) provides some useful bounds on the Gamma operators introduced in Definition \ref{d:gamma}.

\begin{prop}[Estimates on Gamma operators]\label{p:gammaest} For each integer $q\geq 2$ there exists positive constants $c_0,c_1,c_2$ (only depending on $q$)
such that, for all $F=I_q(f)$ with $f\in\HH^{\odot q}$ and $E[F^2]=1$,
one has
\begin{eqnarray}
E\left[\left(\Gamma_2(F)-\frac12\kappa_3(F)\right)^2\right]^{1/2}& \leq&  c_0\,  \kappa_4(F)^{\frac34};\notag\\
E[|\Gamma_3(F)|] &\leq&  c_1\, \kappa_4(F);\notag\\
E[|\Gamma_4(F)|] &\leq&  c_2\, \kappa_4(F)^{\frac54} \label{E1}.
\end{eqnarray}
\end{prop}
{\it Proof}. See \cite[Proposition 4.3]{BBNP}.\qed

\section{Proof of Theorem \ref{t:1}}

Since $ E[F_n^2] = 1$, the fact that $F_n$ converges in distribution to $N\sim \mathscr{N}(0,1)$ if and only $E[F_n^4] \to 3$ is a direct consequence of the main result of \cite{NunuGio}. Whenever $E[F_n^4] \to 3$, one has that the collection of random variables $\{F_n^3 : n\geq 1\}$ is uniformly integrable, in such a way that, necessarily, $E[F_n^3] \to E[N^3] = 0$. The upper and lower bounds appearing in formula \eqref{e:mainest} will be proved separately. 
\smallskip

\noindent ({\it Upper bound} ) We first establish some preliminary estimates concerning a general centered random variable $F\in \mathbb{D}^{1,2}$. We have
\begin{eqnarray*}
\frac12 d_{TV}(F,N)&\leq& E\big[\big|E[1-\langle DF,-DL^{-1}F\rangle|F]\big|\big]  \\
&=& E\Big[E[1-\langle DF,-DL^{-1}F\rangle|F] \times {\rm sign}\{E[1-\langle DF,-DL^{-1}F\rangle|F] \} \Big]
\\
&\leq& \sup_{g\in C^1_c : \|g\|_\infty\leq 1}\! \! \! \! E[g(F)(1-\langle DF,-DL^{-1}F\rangle)], 
\end{eqnarray*}
where $C^1_c$ denotes the class of all continuously differentiable functions with compact support, and we have implicitly applied Lusin's Theorem. Fix a test function $g\in C_c^1$ bounded by 1, and consider a random variable $N\sim \mathscr{N}(0,1)$ independent of $F$. By virtue of Proposition \ref{p:stein}, writing $\varphi = f_g$ for the solution of the Stein's equation associated with $g$ (see \eqref{e:j}), one has that $\|\varphi\|_\infty\leq 2\sqrt{\frac{2}{\pi}}$ and $\|\varphi'\|_\infty\leq 4$. Exploiting independence together with the fact that $E[1-\langle DF,-DL^{-1}F\rangle] = 1-{\rm Var}(F)=0$, we deduce that
\begin{eqnarray*}
&&E[g(F)(1-\langle DF,-DL^{-1}F\rangle)]\\
&=& E[(g(F)-E[g(N)])(1-\langle DF,-DL^{-1}F\rangle)]\\
&=& E[(\varphi'(F)-F\varphi(F))(1-\langle DF,-DL^{-1}F\rangle)]\\
&=&E[\varphi'(F)(1-\langle DF,-DL^{-1}F\rangle)]
-E[\varphi'(F)\langle DF,-DL^{-1}F\rangle(1-\langle DF,-DL^{-1}F\rangle)\\
&&
+E[\varphi(F)\Gamma_2(F)]\\
&=&E[\varphi'(F)(1-\langle DF,-DL^{-1}F\rangle)^2]
+E[\varphi(F)\Gamma_2(F)],
\end{eqnarray*}
where we have used several times the integration by parts formula \eqref{e:ipp}. In order to properly assess the term $E[\varphi(F)\Gamma_2(F)]$, we shall consider the function $\psi = f_\varphi$, corresponding to the solution of the Stein's equation associated with $\varphi$. Using again Proposition \ref{p:stein}, we deduce the estimates $\|\psi\|_\infty\leq \frac{8}{\pi}$ and $\|\psi'\|_\infty\leq 16$, and moreover
\begin{eqnarray*}
&&E[\varphi(F)\Gamma_2(F)]-  \frac12
E[\varphi(N)]\kappa_3(F)\\
&=& E[(\varphi(F)-E[\varphi(N)])\Gamma_2(F)] 
= E[(\psi'(F)-F\psi(F))\Gamma_2(F)]\\
&=&E[\psi'(F)\Gamma_2(F)]
-E[\psi'(F)\langle DF,-DL^{-1}F\rangle\Gamma_2(F)]
+E[\psi(F)\Gamma_3(F)]\\
&=&E[\psi'(F)(1-\langle DF,-DL^{-1}F\rangle)\Gamma_2(F)]
+E[\psi(F)\Gamma_3(F)],
\end{eqnarray*}
where we used once again independence and integration by parts. Combining the previous bounds, one infers that
\begin{eqnarray*}
d_{TV}(F,N)
&\leq& \sqrt{\frac{2}{\pi}}|\kappa_3(F)|+ 2\sqrt{\frac{2}{\pi}} E[(1-\langle DF,-DL^{-1}F\rangle)^2]\\
&&+
16\sqrt{E[(1-\langle DF,-DL^{-1}F\rangle)^2]}\sqrt{E[\Gamma_2(F)^2]}
+\frac{8}{\pi}E[|\Gamma_3(F)|].
\end{eqnarray*}
To conclude, let us consider a sequence $F_n=I_q(f_n)$, $n\geq 1$, living in the $q$th Wiener chaos of $X$ and such that each $F_n$ has variance 1. Assume that $F_n$ converges in distribution to $N$. Then, for $n$ large enough one has that $|\kappa_3(F_n)|\leq 1$ and $\kappa_4(F_n)\leq 1$.
Using Proposition \ref{p:gammaest}, we immediately deduce that, for some universal constant $c_q>0$ (depending only on $q$),
\begin{eqnarray*}
d_{TV}(F_n,N) &\leq& c_q\left(|\kappa_3(F_n)|+ \kappa_4(F_n) + \sqrt{\kappa_4(F_n)}
\sqrt{\kappa_3(F_n)^2 + \kappa_4(F_n)^{3/2}}\right)\\
&\leq& C \,\max(|\kappa_3(F_n)|, \kappa_4(F_n)),
\end{eqnarray*}
where $C$ is the constant appearing in the statement.

\smallskip

\noindent ({\it Lower bound}) According to the representation \eqref{e:tv2}, the distance $d_{TV}(F_n,N)$ is bounded from below by the quantity
\begin{equation}\label{maw}
\frac12\max\big\{\big|E[\cos(F_n)]-E[\cos(N)]\big|,\big|E[\sin(F_n)]-E[\sin(N)]\big|\big\}.
\end{equation}
Combining \cite[Corollary 3.12]{BBNP} with Proposition \ref{p:gammaest}, one can write
\begin{eqnarray*}
\left|E[\sin(F_n)]-E[\sin(N)]-\frac12E[f''_{\sin}(F_n)]\kappa_3(F_n)
-\frac16E[f'''_{\sin}(F_n)]\kappa_4(F_n)\right|&\leq& 2 E|\Gamma_4(F_n)|\\
&\leq& C\,\kappa_4(F_n)^{\frac54}.
\end{eqnarray*}
Here, $C$ denotes a positive constant which is independent of $n$ and whose value can change from line to line, whereas
 $f_{\sin}$ stands for the solution of the Stein's equation associated with the sine function, as given in \eqref{e:j}. From \cite[formula (5.2)]{BBNP}, one has that
 $E[f''_{\sin}(N)]=-\frac13 E[\sin(N)H_3(N)]=\frac{1}{\sqrt{e}}$.
 Similarly, $E[f'''_{\sin}(N)]=-\frac14 E[\sin(N)H_4(N)]=0$.
Moreover, from \cite[Theorem 1.1]{daly} it comes that $f''_{\sin}$ and
 $f'''_{\sin}$ are both bounded by 2.
 Finally, from (\ref{NPfirst}) one has that 
 $d_{TV}(F_n,N)\leq C\sqrt{\kappa_4(F_n)}$.
Combining all these facts leads to
 \begin{eqnarray*}
 \frac14\big|E[f''_{\sin}(F_n)]-E[f''_{\sin}(N)\big|
\leq d_{TV}(F_n,N)&\leq& C\sqrt{\kappa_4(F_n)};\\
 \frac14\big|E[f'''_{\sin}(F_n)]-E[f'''_{\sin}(N)\big|
\leq d_{TV}(F_n,N)&\leq& C\sqrt{\kappa_4(F_n)}
 \end{eqnarray*}
 so that
 \[
\left| E[\sin(F_n)]-E[\sin(N)] -\frac{1}{2\sqrt{e}}\,\kappa_3(F_n)\right|\leq C\,\max \big\{|\kappa_3(F_n)|,\kappa_4(F_n)\big\}\times \kappa_4(F_n)^{\frac14}.
 \]
 Similarly,
 one shows that
  \[
\left| E[\cos(F_n)]-E[\cos(N)] +\frac{1}{4\sqrt{e}}\,\kappa_4(F_n)\right|\leq C\,\max \big\{|\kappa_3(F_n)|,\kappa_4(F_n)\big\}\times \kappa_4(F_n)^{\frac14}.
 \]
 As a consequence, exploiting the lower bound (\ref{maw}) we deduce that
 \[
 d_{TV}(F_n,N)\geq \left(\frac{1}{4\sqrt{e}}-C\kappa_4(F_n)^{\frac14}\right)\max \big(|\kappa_3(F_n)|,\kappa_4(F_n)\big),
 \]
 so that the proof of the theorem is concluded.

\section{Application to the Hermite variations of the discrete-time fractional Brownian motion}

We now discuss an application of Theorem \ref{t:1} to non-linear functionals of a fractional Brownian motion. Consider a fractional Brownian motion with Hurst index $H\in(0,1)$.
We recall that $B_H=\{B_H(t) : t\in \R\}$ is a centered
Gaussian process with continuous paths such that
\[
E[B_H(t)B_H(s)] = \frac12\big(|t|^{2H}+|s|^{2H}-|t-s|^{2H}\big), \quad s,t\in\R.
\]
The process $B_H$ is self-similar with stationary increments. We refer the reader to Nourdin \cite{NouBookFBM} for a self-contained introduction to its main properties.

\medskip

In this section, we shall rather work with the so-called {\it fractional Gaussian noise} associated with $B_H$, which is the Gaussian sequence given by
\begin{equation}\label{xk}
X_k=B_{H}(k+1)-B_H(k),\quad k\in\Z.
\end{equation}
Note that the family $\{X_k : k\in \Z\}$ constitutes a centered stationary Gaussian family with covariance
\[
\rho(k)=E[X_rX_{r+k}]=\frac{1}{2}\left(|k+1|^{2H}-2|k|^{2H}+|k-1|^{2H}\right),\quad r,k\in\Z.
\]
Also, it is readily checked that $\rho(k)$ behaves asymptotically as
$\rho(k)\sim H(2H-1)|k|^{2H-2}$ as $|k|\to\infty$.

\medskip

Now fix an integer $q\geq 2$, consider the $q$th Hermite polynomial $H_q$, and set
\[
F_n := \frac{1}{\sqrt{n\,v_n}}\sum_{k=0}^{n-1} H_q(X_k),
\]
with $v_n>0$ chosen so that $E[F_n^2]=1$. An important problem in modern Gaussian analysis is to characterise those values of $H,q$ such that the sequence $F_n$ verifies a CLT. The following statement is the celebrated {\it Breuer-Major CLT}, first proved in \cite{breuermajor} (see \cite[Chapter 7]{NP11} for a modern proof and for an overview of its many ramifications).

\begin{thm}[Breuer-Major CLT \cite{breuermajor}] Let the previous notation and assumptions prevail. Then, one has that
\begin{equation}\label{cvlaw}
F_n\overset{\rm Law}{\longrightarrow}N \sim\mathscr N(0,1),\quad\mbox{ as $n\to\infty$,}
\end{equation}
if and only if $H$ belongs to the interval $(0,1-\frac{1}{2q}]$. 
\end{thm}

Deducing explicit estimates on the speed of convergence in the CLT \eqref{cvlaw} is a difficult problem, that has generated a large amount of research (see \cite{NouBookFBM, NP11} for an overview of the available literature, as well as \cite{BBL, BBNP, NNT} for recent developments). Using our Theorem \ref{t:1} together with the forthcoming Proposition \ref{p:ppp} allows one to deduce exact rates of convergence in total variation for every value of $q,H$ such that $H \in (0,1-\frac{1}{2q})$. We adopt the following convention for non-negative sequences $(u_n)$ and $(v_n)$: we write $v_n\propto u_n$ to indicate that $0<\liminf_{n\to\infty}v_n/u_n\leq \limsup_{n\to\infty}v_n/u_n<\infty$.

\begin{prop}[See \cite{BBNP}]\label{p:ppp}
Let the above notation and assumptions prevail. 

\begin{enumerate}

\item If $q$ is odd, then $\kappa_3(F_n) = E[F_n^3]=0$ for every $n$.

\item For every even integer $q\geq 2$,
\[
\kappa_3(F_n)\propto
\left\{\begin{array}{lll}
n^{-\frac 12} \;&\mbox{if } &0<H<1-\frac2{3q} \\
n^{-\frac 12}\log^2 n \;&\mbox{if } &H=1-\frac2{3q} \\
n^{\frac32-3q+3qH}\;&\mbox{if }& 1-\frac2{3q}<H<1-\frac{1}{2q} \end{array}\right..
\]

\item For $q\in\{2,3\}$, one has that
\begin{equation}\label{rate-quatre}
\kappa_4(F_n)\propto
\left\{\begin{array}{lll}n^{-1}\;&\mbox{if } & 0<H  <1-\frac{3}{4q} \\
n^{-1}\log^3 n \;&\mbox{if } &H=1-\frac{3}{4q} \\
n^{4qH-4q+2}\;&\mbox{if }& 1-\frac{3}{4q}<H<1-\frac{1}{2q} \end{array}\right.\;.
\end{equation}

\item For every integer $q>3$,
\begin{equation}\label{rate-qBFM}
\kappa_4(F_n)\propto
\left\{\begin{array}{lll}
n^{-1}\;&\mbox{if }& 0<H <\frac34\\
n^{-1}\log (n)\;&\mbox{if }& H =\frac34 \\
n^{4H-4}\;&\mbox{if } &\frac34<H<1-\frac1{2q-2}\\
n^{4H-4}\log^2 n\;&\mbox{if } & H  = 1-\frac 1{2q-2}\\
n^{4qH-4q+2}\;&\mbox{if } &  1-\frac 1{2q-2}<H<1-\frac{1}{2q}
\end{array}\right.\;.
\end{equation}
\end{enumerate}
\end{prop}

Considering for example the cases $q=2$ and $q=3$ yields the following exact asymptotics, which are outside the scope of any other available technique.

\begin{prop} \label{sharp:k3}

\begin{enumerate}
\item If $q=2$, then
\[
d_{TV}(F_n,N)\propto
\left\{\begin{array}{lll}
n^{-\frac 12} \;&\mbox{if } &0<H<\frac2{3} \\
n^{-\frac 12}\log^2 n \;&\mbox{if } &H=\frac2{3} \\
n^{6H-\frac92}\;&\mbox{if }& \frac2{3}<H<\frac{3}{4} \end{array}\right..
\]
\item If $q=3$, then
\[
d_{TV}(F_n,N)\propto
\left\{\begin{array}{lll}
n^{-1}\;&\mbox{if } & 0<H  <\frac{3}{4} \\
n^{-1}\log^3 n \;&\mbox{if } &H=\frac{3}{4} \\
n^{12H-10}\;&\mbox{if }& \frac{3}{4}<H<\frac{5}{6} 
\end{array}\right.\;.
\]
\end{enumerate}
\end{prop}

\bibliographystyle{abbrv}

\begin{thebibliography}{}

\end{thebibliography}


\begin{thebibliography}{10}


\bibitem{BBL} H. Bierm\'e, A. Bonami and J.R. Le\'on (2011). Central limit theorems and quadratic variations in terms of spectral density. {\it Electron. J. Probab.} {\bf 16}(13), 362-395.


\bibitem{BBNP} H. Bierm\'e, A. Bonami, I. Nourdin and G. Peccati (2012). Optimal Berry-Esseen rates on the Wiener space: the barrier of third and fourth cumulants. {\it ALEA} {\bf 9}(2), 473-500.


\bibitem{breuermajor}
P. Breuer and P. Major (1983). Central limit theorems for non-linear functionals of Gaussian fields. {\it J. Mult. Anal.} {\bf 13}, 425-441.

\bibitem{campese} S. Campese (2013). Optimal convergence rates and one-term Edgeworth
expansions for multidimensional functionals of Gaussian fields. In preparation.



\bibitem{ChGoSh11}
L.~H.~Y. Chen, L.~Goldstein, and Q.-M. Shao (2011).
\newblock {\em Normal approximation by {S}tein's method}.
\newblock Probability and its Applications (New York). Springer, Heidelberg.


\bibitem{daly}
F.A. Daly (2008).
Upper bounds for Stein-type operators.
{\it Electron. J. Probab.} {\bf 13}, paper 20, 566-587.


\bibitem{dudley_book}
R.~M. Dudley (1989).
\newblock {\em Real analysis and probability}.
\newblock The Wadsworth \& Brooks/Cole Mathematics Series. Wadsworth \&
  Brooks/Cole Advanced Books \& Software, Pacific Grove, CA.



\bibitem{MP11}
D.~Marinucci and G.~Peccati (2011).
\newblock {\em Random fields on the sphere : Representation, limit theorems and
  cosmological applications}, volume 389 of {\em London Mathematical Society
  Lecture Note Series}.
\newblock Cambridge University Press, Cambridge.

\bibitem{nourdinln}I. Nourdin (2013). Lectures on Gaussian approximations with Malliavin calculus. To appear in: {\it S\'eminaire de probabilit\'es}.

\bibitem{NouBookFBM}
I.~Nourdin (2012).
\newblock {\em Selected Aspects of Fractional Brownian Motion}.
\newblock Springer-Verlag.



\bibitem{NNT} I. Nourdin, D. Nualart and C.A. Tudor (2010). Central and non-central limit theorems for weighted power variations of fractional Brownian motion. {\it Ann. I.H.P.} {\bf 46}(4), 1055-1079.

\bibitem{np-ptrf}
\rm I. Nourdin and G. Peccati (2009):
\rm Stein's method on Wiener chaos.
\it Probab. Theory Relat. Fields \rm {\bf 145}, no. 1, 75-118.

\bibitem{np-aop}
\rm I. Nourdin and G. Peccati (2009):
\rm Stein's method and exact Berry-Esseen asymptotics for functionals of Gaussian fields.
\it Ann. Probab. \rm {\bf 37}, no. 6, 2231-2261.


\bibitem{np-jfa}
\rm I. Nourdin and G. Peccati (2010).
\rm Cumulants on the Wiener space.
\it J. Funct. Anal. \rm {\bf 258}, 3775-3791.


\bibitem{NP11}
I.~Nourdin and G.~Peccati (2012).
\newblock {\em Normal approximations with Malliavin calculus : from Stein's
  method to universality}.
\newblock Cambridge Tracts in Mathematics. Cambridge University Press.

\bibitem{NPRaop}
I.~Nourdin, G.~Peccati, and G.~Reinert (2010).
\newblock Invariance principles for homogeneous sums: universality of
  {G}aussian {W}iener chaos.
\newblock {\em Ann. Probab.} {\bf 38}(5):1947--1985.

\bibitem{NoPeSw} I. Nourdin, G. Peccati and Y. Swan (2013). Entropy and the fourth moment phenomenon. Preprint.


\bibitem{NV-ejp}
I. Nourdin and F.~G. Viens (2009).
\newblock Density formula and concentration inequalities with Malliavin calculus. 	
\newblock {\em Electron. J. Probab.} {\bf 14}, 2287--2309. 

\bibitem{Nu06}
D.~Nualart (2006).
\newblock {\em The {M}alliavin calculus and related topics}.
\newblock Probability and its Applications (New York). Springer-Verlag, Berlin,
  second edition.


\bibitem{NunuGio}
D.~Nualart and G.~Peccati (2005).
\newblock Central limit theorems for sequences of multiple stochastic
  integrals.
\newblock {\em Ann. Probab.} {\bf 33}(1), 177--193.

\bibitem{viensSPA}
F.~G. Viens (2009).
\newblock Stein's lemma, {M}alliavin calculus, and tail bounds, with
  application to polymer fluctuation exponent.
\newblock {\em Stochastic Process. Appl.} {\bf 119}(10), 3671--3698.

\end{thebibliography}

\end{document}